\documentclass[12pt]{article}
\usepackage{xcolor,graphicx} 
\usepackage{amsmath,amsthm,amssymb}
\newtheorem{theorem}{Theorem}[section]

\newtheorem{lemma}[theorem]{Lemma}
\newtheorem{proposition}[theorem]{Proposition}

\newtheorem{claim}[theorem]{Claim}

\newtheorem{quest}[theorem]{Question}

\def\beq{\begin{equation}}\def\eeq{\end{equation}}
\def\beqn{\begin{eqnarray}}\def\eeqn{\end{eqnarray}}

\def\qed{\ifhmode\unskip\nobreak\fi\quad\ifmmode\Box\else$\Box$\fi}

\title{Clique covers of complete graphs and piercing multitrack intervals}
\author{J\'anos Bar\'at\thanks{Alfr\'ed R\'enyi Institute of Mathematics, Budapest, P.O. Box 127, Budapest, Hungary, H-1364.
\texttt{barat.janos@renyi.hu}, \texttt{gyarfas.andras@renyi.hu}, \texttt{sarkozy.gabor@renyi.hu} }\footnotemark[1]
\thanks{University of Pannonia, Research supported in part by NKFIH, Grant Number K131529 and ERC Advanced Grant ``GeoScape'' No. 882971.} \and Andr\'as Gy\'arf\'as\footnotemark[1] \thanks{Research supported in part by
NKFIH Grant No. K132696.} \and G\'{a}bor N. S\'ark\"ozy\footnotemark[1]
\thanks{Computer Science Department, Worcester Polytechnic Institute, Worcester, MA.} \thanks{Research supported in part by
NKFIH Grants No. K132696, K117879.}}

\begin{document}
\maketitle
\date{}
\begin{abstract} Assume that $R_1,R_2,\dots,R_t$ are disjoint parallel lines in the plane. A $t$-interval (or $t$-track interval) is a set that can be written as the union of $t$ closed intervals, each on a different line. It is known that pairwise intersecting $2$-intervals can be pierced by two points, one from each line.
However, it is not true that every set of pairwise intersecting $3$-intervals can be pierced by three points, one from each line.
For $k\ge 3$, Kaiser and Rabinovich asked whether $k$-wise intersecting $t$-intervals can be pierced by $t$ points, one from each line. Our main result provides a positive answer in an asymptotic sense: in any set $S_1,\dots,S_n$ of $k$-wise intersecting $t$-intervals, at least $\frac{k-1}{k+1}n$ can be pierced by $t$ points, one from each line. We prove this in a more general form, replacing intervals by subtrees of a tree.
This leads to questions and results on covering vertices of edge-colored complete graphs by vertices of monochromatic cliques having distinct colors, where the colorings are chordal, or more generally induced $C_4$-free graphs. For instance, we show that if the edges of a complete graph $K_n$ are colored with red or blue so that both color classes are induced $C_4$-free, then at least ${4n\over 5}$ vertices can be covered by a red and a blue clique, and this is best possible. We conclude by pointing to new Ramsey-type problems emerging from these restricted colorings.

 \end{abstract}

\section{Introduction, results}

\subsection{$t$-intervals}
Assume that $R_1,R_2,\dots,R_t$ are disjoint parallel lines in the plane. A $t$-interval (or $t$-track interval) is a set that can be written as the union of $t$ closed intervals, $I_1,I_2,\dots,I_t$ with $I_j\subset R_j$. Gallai (in connection with a problem of piercing directed cycles in a digraph) asked how many points are needed to pierce pairwise intersecting $2$-track intervals. Here we use  {\em piercing} as an equivalent notion of {\em cover}: a family $H$ of subsets of a set $V$ is pierced (or covered) by $T\subset V$ if every set in $H$ has a non-empty intersection with $T$.

Gy\'arf\'as and Lehel proved in \cite{GYL} that the piercing number of pairwise intersecting $t$-track intervals is bounded by a function of $t$ and for $t=2,3$ the smallest value of the piercing number is $2$ and $4$, respectively. Initiated by a breakthrough result of Tardos \cite{T}, topological methods and their simplifications resulted in sharp and approximate results about the smallest possible piercing number of $t$-intervals where the condition ``pairwise intersecting'' was replaced by ``at most $k$ pairwise disjoint''. For a survey of these results see Chapter 30 in \cite{B}, Chapter 3 in \cite{MZB}.

The motivation of this paper is an interesting special piercing of $t$-intervals introduced by Kaiser and Rabinovich \cite{KR}: the piercing
set is required to contain {\em at most one point} from every $R_i$. They call a family of $t$-track intervals {\em weakly intersecting} if there exists such a piercing set.
A set system is {\em $k$-wise intersecting} or {\em $k$-intersecting}, if any $k$ sets have a non-empty intersection.
Replacing pairwise intersection in the previous problem by $k$-wise intersection, the following was asked in \cite{KR}.

\begin{quest}[Kaiser, Rabinovich]\label{krquestion} Assume that $k\ge 3$ and $H$ is a family of $k$-wise intersecting $t$-intervals. Is $H$ weakly intersecting?
 \end{quest}

For $k\geq \lfloor \log_2(t+1)\rfloor+1$, Kaiser and Rabinovich proved (Theorem 4.5 in \cite{KR}) that $k$-wise intersecting $t$-intervals are weakly intersecting. In particular, $3$-wise intersecting $t$-intervals are weakly intersecting if $t\le 6$.

Here we focus on a natural refinement towards Question \ref{krquestion}: given a set $H$ of $k$-wise intersecting $t$-intervals, how large subset $F$ of $H$ is weakly intersecting? Our main result shows that (independently of $t$), $|F|$ approaches to $|H|$ as $k$ grows, thus provides a positive answer to Question \ref{krquestion} in an asymptotic sense.

\begin{theorem}\label{intcor} Let $t$ and $k\ge 2$ be positive integers. If $H$ is a set of $k$-wise intersecting $t$-intervals, then at least $\frac{(k-1)|H|}{k+1}$ are weakly intersecting.
\end{theorem}

We prove Theorem \ref{intcor} in a more general form, instead of intervals, for subtrees of a tree (Theorem \ref{lower}).
We also show (see the discussion at the end of the proof of Theorem \ref{lower}) that a similar statement holds if we replace $k$-wise intersecting with a different
strengthening of pairwise intersecting. Namely, with {\em pairwise $k$-fold intersection}, that is, every pair of $t$-intervals
intersect on at least $k$ of the lines.
For $k\geq \lfloor \frac{t+1}{2} \rfloor$, Kaiser and Rabinovich proved (Theorem 4.6 in \cite{KR}) that  pairwise $k$-fold intersecting $t$-intervals are weakly intersecting.
We focus on $k$-wise intersection in this paper. However, the corresponding problems are interesting for pairwise $k$-fold intersection as well.

Theorem \ref{intcor} with $k=2$ shows that there is $F\subset H, |F|\ge {|H|\over 3}$ such that $F$ is weakly intersecting. Our next result complements this with an upper bound on $|F|$.

\begin{theorem}\label{onefourth}  For any $t\geq 2, n\ge 4t-5$, there exists a set $H$ of pairwise intersecting $t$-intervals such that $|H|=n$ and for any weakly intersecting subset $F\subset H$, we have $|F|\le {3(t-1)\over 4t-5}n$.
\end{theorem}

\begin{figure}[h!]
\begin{center}
 \includegraphics[width=0.7\linewidth]{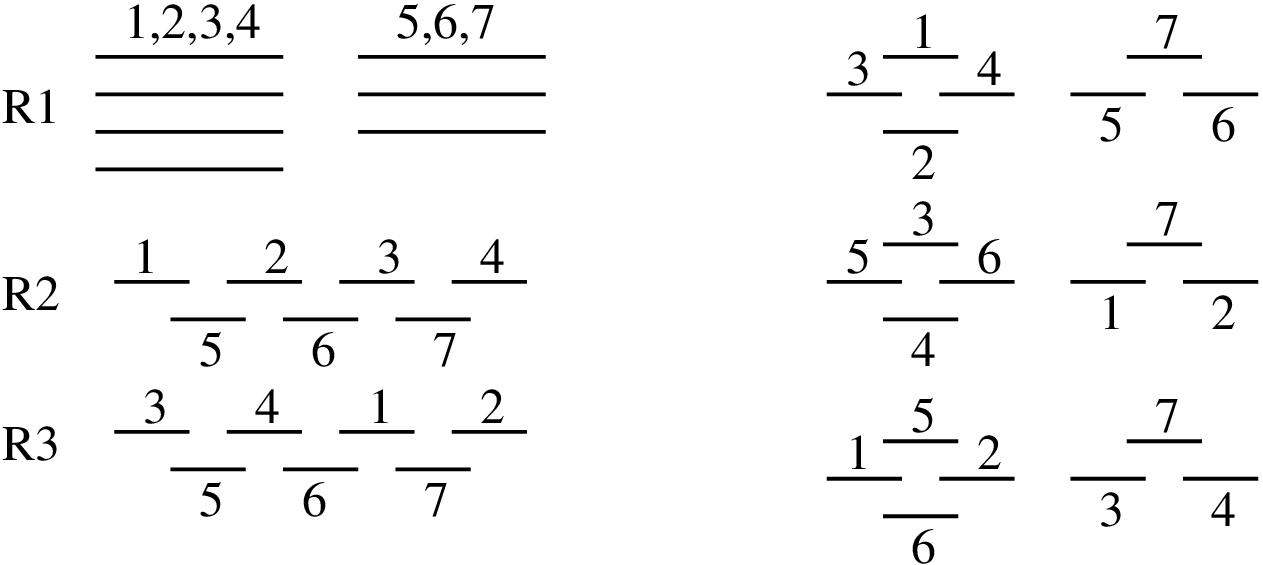}
 \caption{\label{c1c2}The constructions of \cite{GYL} and \cite{KR}.}
\end{center}
\end{figure}

Pairwise intersecting $3$-intervals are not necessarily weakly intersecting, constructions from \cite{GYL} and \cite{KR} are shown in Figure~\ref{c1c2}.  Note that the first example can be pierced by two points on  $R_1$ (but not with three points from distinct lines).

\subsection{Strong monochromatic clique covers}
The geometric problems mentioned so far can be asked in a more general framework. Consider a complete graph $K_n$, whose edge set is the union of $t$ graphs, $G_1,G_2,\dots,G_t$. We consider this as an edge-coloring of $K_n$, where edges can have several colors, subsets of $[t]$. A $(t,k)$-coloring of $K_n$ is a (multi)coloring with $t$ colors, where each set of $k$ vertices spans a monochromatic complete subgraph in some color. We often use the shorter term {\em clique} for complete graphs.

A covering of the vertices of $K_n$ with vertices of monochromatic cliques is called a {\em strong cover} if the colors of the covering cliques are all different. If a $t$-coloring of $K_n$ has a strong cover, then we define $\theta(K_n)$ as the minimum number of cliques needed for a strong cover. (Note that a strong cover can be reduced to a strong partition.)

One can easily see that a set of $k$-wise intersecting $t$-intervals is equivalent to a $(t,k)$-coloring, where each $G_i$ is an {\em interval graph} i.e. can be represented as the intersection graph of a set of intervals on a line. We call such a coloring an {\em interval $(t,k)$-coloring}. Also, ``weakly intersecting'' corresponds to ``strong covering'' by the one-dimensional Helly theorem, a clique monochromatic in color $i$ corresponds to a piercing point on $R_i$.
A well-known extension of interval graphs is the family of {\em chordal (or triangulated)} graphs, that do not contain cycles of length at least four as an induced subgraph. Chordal graphs are interesting from a geometric point of view as well, they can be characterized as intersection graphs of subtrees of a tree \cite{GAV}. A $(t,k)$-coloring is {\em chordal} if each $G_i$ is a chordal graph. Thus a strong cover corresponds to weakly intersecting $t$-subtrees. Properties of chordal graphs can be found in  \cite{GO}. A further extension is to consider {\em $C_4$-free graphs}, where only induced cycles of length four are forbidden. In what follows, we use $C_4$-free, $C_5$-free in the induced sense.  A simple but useful property of interval/chordal/$C_4$-free colorings is that they are closed for clique substitution. A {\em clique substitution} into a $(t,k)$-colored $K_n$ is defined as follows. Replace a vertex $v$ of $K_n$ by a clique $K$ colored with all the $t$ colors. All edges between $K$ and $u\in V(K_n)$ are colored with the colors on the edge $(u,v)\in E(K_n)$. In case of $V(K)=\{v,w\}$ we say that $w$ is a {\em replica} of $v$.

\begin{proposition}\label{subst} Assume that $\theta(K_n){=}s$ in an interval/chordal/$C_4$-free $(t,k)$-coloring of $K_n$ and $K_m$ is obtained by a clique substitution into $K_n$. The obtained coloring is an interval/chordal/$C_4$-free $(t,k)$-coloring and $\theta(K_m)=s$.
\end{proposition}

We prove Theorem \ref{intcor} in the following stronger form.

\begin{theorem} \label{lower}
Let $t,k\geq 2$ and consider a chordal $(t,k)$-coloring of the complete graph $K_n$. There is a strong cover on at least $\frac{k-1}{k+1}n$ vertices of $K_n$.
\end{theorem}

For chordal $(2,2)$-colorings of $K_n$ there is a strong cover for all vertices \cite{GYL}. This holds for chordal $(3,3)$-colorings of $K_n$ as well. (For interval $(3,3)$-colorings Theorem \ref{33} is proved in \cite{KR}.)

\begin{theorem}\label{33}   Assume that we have a chordal $(3,3)$-coloring of $K_n$. Then $\theta(K_n)\le 3$ and {\rm Figure~\ref{abra3,3}} shows that equality can occur for interval graphs.
\end{theorem}

\begin{figure}[h!]
\begin{center}
 \includegraphics[width=0.3\linewidth]{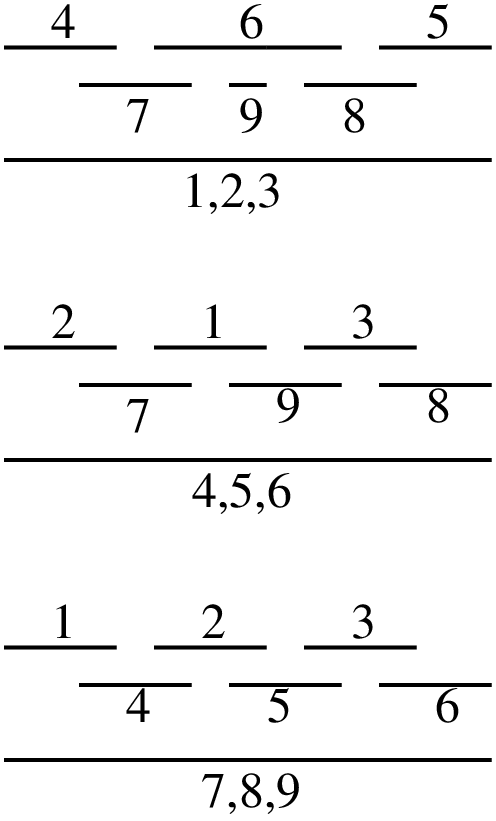}
   \caption{\label{abra3,3} Nine $3$-wise intersecting $3$-intervals without piercing by two points.}
\end{center}
\end{figure}

\begin{theorem}\label{cortt} Assume that we have a chordal $(t,t)$-coloring of $K_n$. For even $t\ge 2$, $\theta(K_n)=2$, and for odd $t\ge 3$, $\theta(K_n)\le 3$.
\end{theorem}

For $C_4$-free $(2,2)$-colorings, we can determine the best strong cover of $K_n$.

\begin{theorem}\label{fourover5}  In every $C_4$-free $(2,2)$-coloring of $K_n$ there is a strong cover on some $K_m\subset K_n$ with $m\ge \lceil{4n\over 5}\rceil$. This is best possible for every $n\ge 5$.
\end{theorem}




\section{Proofs}
\noindent
{\bf Proof of Theorem \ref{onefourth}.} We present the construction (a generalization of Figure~\ref{c1c2} left part) in the framework of interval $(t,2)$-colorings. Set $n=4t-5$, and define the following $t$-coloring of $K_n$.
Edges of color $1$ form two vertex-disjoint cliques,  $A=V(K_{2t-2})$ and $B=V(K_{2t-3})$, this is clearly an interval graph. The other $t-1$ colors are defined by decomposing the complete bipartite graph $[A,B]$ into $t-1$ Hamilton paths (which are also interval graphs). The easiest way to get this partition is to extend $B$ with one vertex and use the well-known fact that $K_{2t-2,2t-2}$ can be partitioned into Hamilton cycles (an early reference is \cite{AL}).  Removing one vertex from $B$, we get the required partition into Hamilton paths.

Observe that a partition of $V(K_n)$ into cliques of different colors can use only one of the two cliques in the first color. However, the vertices of the unused clique are independent in all other colors. Therefore, we can cover at most $2t-2+t-1=3(t-1)$ vertices by $t$ cliques of distinct colors.
At least $t-2$ vertices are uncovered, giving the ratio ${3(t-1)\over 4t-5}$ for the covered vertices. Applying Proposition \ref{subst}, we can extend this example to an interval $(t,2)$-coloring of $K_n$, for any $n\ge 4t-5$.   \qed

\medskip

\noindent
{\bf Proof of Proposition \ref{subst}.}  To see that substituting a clique $K$ to $v\in V(K_n)$ preserves $(t,k)$-coloring, consider any $S\subset V(K_m)$ with $|S|=k$. Set $M=S\cap K$. If $|M|\le 1$, then $S$ spans a monochromatic clique by the assumption that $K_n$ has a $(t,k)$-coloring. Otherwise, for any $w\in M$, $S'=M\setminus \{w\}$ extends the monochromatic clique $S\setminus S'$ by the definition of the substitution.

Substituting a clique $K$ can create a chordless cycle $C$ of length at least $4$ only if $|V(C)\cap V(K)|=2$, and the intersection is an edge $e=(v_i,v_{i+1})$  of $C$. Indeed, if $|V(C)\cap V(K)|\ge 3$, then $C$ cannot be chordless.  If $|V(C)\cap V(K)|\le 1$, then $C$ was present before the substitution. However, this implies that the next vertex $v_{i+2}$ on $C$ is not in $K$, thus adjacent to $v_i$ by the rule of substitution and creates a chord in $C$. Thus clique substitution preserves the property of being interval/chordal/$C_4$-free.

It is obvious that $\theta$ cannot decrease by clique substitution. But it cannot increase either, since the clique $C$ covering the vertex $v$ in an optimal strong cover of $K_n$ is extended by the substituted $K$ into a clique of the same color. \qed

\medskip
\noindent
{\bf Proof of Theorem \ref{lower}.}  The following proof idea is partly inspired by Theorem 3 of \cite{GYHS}.
Consider a chordal $(t,k)$-coloring of the complete graph $K_n$.
Let $G_i$ denote the chordal graph on $n$ vertices induced by color $i$.
Let us consider color~1, and let $\omega_1$ be the clique number of $G_1$.
Since $G_1$ is perfect, the $n$ vertices of $K_n$ can be partitioned into $\omega_1$ independent sets $X_1,\dots, X_{\omega_1}$. Since $G_1$ is chordal, the edges of color $1$ between $X_s$ and $X_z$ form an acyclic subgraph of $K_n$ for any two distinct indices $s,z$. Indeed, otherwise a cycle of minimum length in the bipartite graph $[X_s,X_z]$ would define a chordless cycle in $G_1$. 
Therefore, we get an upper bound on the number of edges of color $1$ induced by any subset of vertices.
For the entire vertex set, we get:

$$|E(G_1)|\le \sum_{s<z} (|X_s|+|X_z|-1)=(\omega_1-1)n-\binom{\omega_1}{2}$$
$$\leq \omega_1 n - \omega_1 = \omega_1 (n-1). $$

Similarly, if we consider an arbitrary subset $R$ of the vertices, then (using the notation $G|_R$ for the subgraph induced by $R$),
$$|E(G_1|_R)|\le \omega_1 (|R|-1).$$

Indeed, if $\omega_1'$ denotes the number of those $X_i$-s that have a non-empty intersection with $R$, we get

$$|E(G_1|_R)|\leq (\omega_1'-1)|R|-\binom{\omega_1'}{2}\leq \omega_1' |R| - \omega_1' = \omega_1' (|R|-1)\leq \omega_1 (|R|-1).$$

In step 1, we select a clique $W_1$ of color $1$ and size $\omega_1$ and remove the vertices of $W_1$ from $K_n$ to get $G^{(2)}$, a $t$-edge-colored complete graph on $n-\omega_1$ vertices. We repeat this step for the remaining colors $2,\dots,t$.

In step 2, we select a clique $W_2$ in $G^{(2)}$ of color 2 and size $\omega_2$ and remove the vertices of $W_2$ from $G^{(2)}$ to get $G^{(3)}$. We also deduce that the number of edges of color 2 on any remaining vertex set $R$ satisfies:
$$|E(G_2|_R)|\le \omega_2 (|R|-1).$$

In step $i$, we select a clique $W_i$ in $G^{(i)}$ of color $i$ and size $\omega_i$ and remove the vertices of $W_i$ from $G^{(i)}$ to get $G^{(i+1)}$. We also deduce that the number of edges of color~$i$ on any remaining vertex set $R$ satisfies:

$$|E(G_i|_R)|\le \omega_i(|R|-1).$$

Suppose we covered $cn$ vertices of $K_n$ after $t$ steps with cliques $W_1,\dots,W_t$ for some constant $c$. Denote the set of uncovered vertices by $T$.
We get an upper bound on $M$, the sum of multiplicities of edges induced in $T$ by some $G_i$, adding the upper bounds for each color.

$$M=\sum_{i=1}^t |E(G_i|_T)|\le \sum_{i=1}^t \omega_i(|T|-1)= \sum_{i=1}^t |W_i|(|T|-1)=cn(|T|-1).$$

On the other hand, we can get a lower bound on $M$ using the following claim.

\begin{claim}\label{fedett}
In $K_n|_T$ every edge has at least $k-1$ colors.
\end{claim}

Indeed, for $k=2$, the claim is trivial, since $K_n$ is a complete graph. For $k\geq 4$, we can prove a slightly stronger claim, namely that {\em every} edge in $K_n$ (not just in $K_n|_T$) has at least $k{-}1$ colors.
Otherwise $K_n$ has a strong cover by at most 3 monochromatic cliques. Indeed, suppose to the contrary that an edge $e=(u,v)$ in $K_n$ has only colors $1, \ldots, k-2$. Consider the subgraph  $G'=G_1\cup \ldots \cup G_{k-2}$ on $V(K_n)$.  We claim that $G'$ is $(k{-}2)$-intersecting (in {\em these} $k-2$ colors). Indeed, consider a set $S\subset V(K_n)$ of $k-2$ vertices. We add $\{u,v\}$ and potentially 1 or 2 more vertices (if $S\cap \{u,v\}\not= \emptyset$) to
get a set $S'$ of $k$ vertices. We had a chordal $(t,k)$-coloring by assumption, so $S'$ is contained in a monochromatic clique. Since this clique contains the edge $e$, this could only be one of the colors $[k-2]$. Thus indeed $S$ is contained in a monochromatic clique in these colors. Furthermore, this implies that every edge of $K_n$ is covered by $G'$.
However, $G'$ is $(k{-}2)$-colored and $(k{-}2)$-intersecting, therefore by  Theorem~\ref{cortt} $V(G')=V(K_n)$ can be covered by at most 3 monochromatic cliques, a much stronger result than what we are looking for.

For $k=3$, suppose to the contrary that edge $e=(u,v)$ has only color $i$ in $K_n|_T$, say $u\in X_s$ and $v\in X_z$ in the decomposition above. Since $u$ is in $T$ (after removing $W_i$), originally $|X_s|\geq 2$ was true.
Add a vertex $u'\in X_s$, $u'\not= u$ to get $S=\{u,u',v\}$.
The graph $K_n$ was 3-intersecting, so $S$ is contained in a monochromatic clique. However, this could only be color $i$, since $e$ only has color $i$. On the other hand, the edge $(u,u')$ cannot have color $i$ since $X_s$ is an independent set in color $i$, a contradiction, finishing the proof of the claim.

Finally, comparing the bounds on $M$, we get

\beq\label{elszam}
(k-1)\binom{(1-c)n}{2}\le M\le cn(|T|-1) = cn ((1-c)n-1).
\eeq

Hence
$$(k-1)(1-c) \leq 2c,$$
or
$$\frac{k-1}{k+1}\leq c.\qed $$
\noindent

As we mentioned in the introduction, we can prove a statement similar to Theorem \ref{intcor} and Theorem \ref{lower} in a slightly different
setting, where $k$-wise intersection is replaced with pairwise $k$-fold intersection. More precisely, we have the following.

\begin{theorem} \label{lower1}
Let $t,k\geq 2$, and consider a chordal $t$-coloring of the complete graph $K_n$, where every edge has at least $k{-}1$ colors. Then there is a strong cover of at least $\frac{k-1}{k+1}n$ vertices of $K_n$.
\end{theorem}

The proof is virtually identical to that of Theorem \ref{lower} above. For the upper bound on $M$ in (\ref{elszam}) we just used the fact that we have a chordal $t$-coloring.
For the lower bound we used Claim \ref{fedett}, which now holds by assumption. The rest of the proof is the same.

\medskip
\noindent
{\bf Proof of Theorem \ref{33}.} Consider a chordal $(3,3)$-coloring of $K_n$. Let $E_j$ denote the set of edges colored with $c_j$ for $j=1,2,3$.
We can either find an edge $e=(u,v)\in E_1$ colored only with $c_1$, or $c_2\cup c_3$ is a chordal $(2,2)$-coloring of $K_n$. In the latter case, $V(K_n)$ can be covered by two monochromatic cliques by Theorem~3 of \cite{GYL}. In the former case, we notice that $T=\{u,v,w\}$ with any $w\notin \{u,v\}$ spans a monochromatic triangle, which must be in color $c_1$.
Therefore, $E_1$ induces a connected graph in $c_1$. If it is a clique, then one monochromatic clique covers $V(K_n)$. Otherwise, we can apply a basic property of chordal graphs: there is a clique cut-set $Q$, i.e. $V(K_n)$ has a non-trivial partition into $A,Q,B$ such that $Q$ spans a clique in  $c_1$ and no edge of the complete bipartite graph $[A,B]$ has color $c_1$. Note that neither $A$ nor $B$ induces an edge $e$ having color $c_1$ only. Otherwise we could extend $e$ to a monochromatic triple intersecting both $A$ and $B$. This must be of color $c_1$, a contradiction since no edge between $A,B$ can be of color $c_1$.
Hence $c_2\cup c_3$ is a chordal 2-coloring of the complete graph on $A\cup B$. Again, by  Theorem~3 of \cite{GYL}, $A\cup B$ can be covered by 2 monochromatic cliques, one in color $c_2$ and one in color $c_3$. Therefore, together with $Q$ (in color $c_1$), $V(K_n)$ has a strong cover with three colors.

Equality is possible for $n=9$, because the interval $(3,3)$-coloring represented in Figure~\ref{abra3,3} cannot be pierced by two points from different lines (in fact two points from the same line cannot pierce them either).    The example can be extended by Proposition \ref{subst} for any $n>9$. \qed

\medskip\noindent

{\bf Proof of Theorem \ref{cortt}.} Assume that $t$ is even. For any pair of colors, there is an edge $e$ of $K_n$ missing both colors, otherwise the chordal coloring of  $K_n$ can be considered as a $(2,2)$-coloring and $\theta(K_n)\le 2$ follows from Theorem 3 in \cite{GYL} and the proof is finished. Thus there are edges $e_1,\dots,e_{{t\over 2}}$, each missing a different pair of colors. However, their union has at most $t$ vertices, which cannot span a monochromatic clique, contradicting the assumption that we have a $(t,t)$-coloring.

For odd $t$, the argument is similar, but here we partition the colors into ${t-3\over 2}$ distinct pairs and one triple. Again, for any pair of colors, there is an edge $e_i$ of $K_n$ missing both colors, otherwise the chordal coloring of  $K_n$ can be considered as a $(2,2)$-coloring and $\theta(K_n)\le 2$ follows from Theorem 3 in \cite{GYL} and the proof is finished.  Also, there exists a triangle $T$ in $K_n$, which is not a monochromatic triangle in any of the three colors of the color-triple, otherwise $\theta(K_n)\le 3$ follows from Theorem \ref{33}, finishing the proof. However,
$$T\cup e_1\cup \dots e_{{t-3\over 2}}$$ covers at most $t$ vertices in $K_n$, which cannot span a monochromatic clique, contradicting the assumption that we have a $(t,t)$-coloring. \qed

\medskip

\noindent
{\bf Proof of Theorem \ref{fourover5}.}  We apply the following theorem of Gy\'arf\'as and Lehel (which strengthens a result of Aharoni, Berger, Chudnovsky and Ziani \cite{ABCZ}).  Let $K_5^*$ denote the $2$-colored $K_5$, where each edge has exactly one color and both color classes form a $C_5$.

\begin{theorem}(\cite{GYL2}))\label{gyl2} Assume that $K_n$ has a $(2,2)$-coloring containing no (induced) $C_4$ and $K_5^*$ subgraphs. Then $\theta(K_n)\le 2$.
\end{theorem}

By Theorem \ref{gyl2}, we may assume that $K_n$ contains a copy of $K_5^*$ (otherwise we have a strong cover on $K_n$). Let $K$  be a maximal subgraph of $K_n$, which can be obtained from $K_5^*$ by substituting cliques with vertex sets  $X_1,\dots,X_5$ to its vertices, i.e. no further vertex substitution into $K$ is a subgraph of $G_n$. Assume w.l.o.g. that the complete bipartite graphs  $[X_i,X_{i+1}]$ ($[X_i,X_{i+2}]$) have only red (only blue) edges in$\pmod{5}$ index arithmetic (and the edges within each $X_i$ have both colors).

\begin{lemma}\label{blowup} The edges between any $w\in V(G_n)\setminus V(K)$  and $V(K)$ have a common color.
\end{lemma}
\medskip

\noindent
{\bf Proof of Lemma \ref{blowup}.} Suppose to the contrary that there exists $w\in V(G_n)\setminus V(K)$ and two vertices, $p,q\in V(K)$ such that $(w,p)$ is red only and $(w,q)$ is blue only.
We may assume $p\in X_1$. If $q\in X_1$, then we can find $q'\in V(K)\setminus X_1$ such that $(w,q')$ is also blue only.
Indeed, otherwise for any $v_2\in X_2, v_5\in X_5$ the edges  $(w,v_2),(w,v_5)$  have color red and $(w,v_2,q,v_5)$ is a red induced $C_4$, a contradiction. Thus we may assume  $q\in V_2$ or $q\in V_3$. Due to the symmetry of the coloring of $K$, we may assume w.l.o.g. that $q\in V_2$.
We claim that $w$ is a replica of $v_5$ in the $K_5=\{p,q,v_3,v_4,v_5\}$, where $v_3\in X_3,v_4\in X_4,v_5\in X_5$ are arbitrary. Note that $(w,v_3)$ is blue only, otherwise $(w,p,q,v_3)$ is an induced red $C_4$. Also, $(w,v_4)$ is red only, otherwise $(w,v_4,p,v_3)$ is an induced blue $C_4$. Moreover, $e=(w,v_5)$ must be colored with both red and blue. Indeed, if $e$ is red only, then $(w,q,v_5,v_3)$ spans a blue induced $C_4$; if $e$ is blue only, then $(w,p,v_5,v_4)$ spans a red induced $C_4$, proving the claim.

Since the choice of $v_3,v_4,v_5$ were arbitrary, the claim implies that all edges from $w$ to $X_3$ are blue only, to $X_4$ red only, to $X_5$ both red and blue. We show that all edges from $w$ to $X_1$ are red only and all edges from $w$ to $X_2$ are blue only. This is true for $p\in X_1$ and $q\in X_2$.
Suppose to the contrary  $v_1\in X_1,v_1\ne p$, and the edge $(w,v_1)$ is also blue.
For any $v_4\in X_4$, $(w,v_1,v_4,q)$ is an induced blue $C_4$, a contradiction. Similarly, if for $v_2\in X_2,v_2\ne q$ the edge $(w,v_2)$ is also red, then for any $v_3\in X_3,v_4\in X_4$, $(w,v_2,v_3,w_4)$ is an induced red $C_4$, a contradiction.

Thus $w$ is a replica of $v_5\in K$ contradicting the assumption that $K$ is maximal, proving the lemma.  \qed

By Lemma \ref{blowup}, $V(K_n)\setminus V(K)$  can be partitioned into $R,B$, where all edges of the complete bipartite graph $[V(K),R]$ are red and all edges of the complete bipartite graph $[V(K),B]$ are blue. Also, any edge $(r_1,r_2)$ in $R$ is red, otherwise for any $v_1\in X_1,v_3\in X_3$, $(r_1,v_1,r_2,v_3)$ is an induced red $C_4$. The same argument can be applied to show that all edges within $B$ are blue. Let $|X_i|=\min\{|X_j|: j\in [5]\}$. Now $X_{i+2}\cup X_{i-2}\cup R$ is a red clique and $X_{i+1}\cup X_{i-1}\cup B$ is a blue clique, and their union covers all but at most $|V_i|\le n/5$  vertices of $K_n$. \qed

\section{Conclusion}
Theorem \ref{intcor}  is a step towards deciding Question \ref{krquestion} motivating our research. Its more general form, Theorem \ref{lower} leads to further questions on strong coverings.
Also, clique covers lead to Ramsey-type questions in interval/chordal/$C_4$-free $(t,2)$-colorings of $K_n$. For these restricted colorings the $t$-color classical Ramsey numbers are linear in $n$. Indeed, chordal and $C_4$-free graphs with $n$ vertices and $c{n\choose 2}$ edges contain cliques of size $c'n$. The best value of $c'$ for chordal graphs is $\alpha(c)=1-\sqrt{1-c}$, see \cite{AK}, \cite{GYHS}. For $C_4$-free graphs the existence of $c'$ is established in \cite{GYHS}, and its best value we know is $\alpha^2(c)$, proved by Holmsen \cite{HO}.

\subsection{Strong covers}

\begin{itemize}
\item (i) We do not know what is the largest subset of weakly intersecting $t$-intervals among $n$ pairwise intersecting $t$-intervals.  Theorems \ref{intcor}, \ref{onefourth} place it between  ${n\over 3}$ and ${3n\over 4}$.

\item (ii) Perhaps $\theta(K_n)\le 2$ for chordal $(t,t)$-colorings of $K_n$ for large enough odd $t$, improving Theorem \ref{cortt}. Figure~\ref{abra3,3} shows that this is false for $t=3$.

\item (iii) For general $t$, $C_4$-free $(t,2)$-colorings of $K_n$ are essentially different from chordal $(t,2)$-colorings, where there are $n/3$ vertices with a strong cover by Theorem \ref{lower}. In a $C_4$-free $(t,2)$-coloring the size of the covered part must depend on $t$ as well. Indeed, Chung and Graham \cite{CG} proved that the $t$-color Ramsey number of $C_4$ is asymptotic to $t^2$.  Thus there exist $C_4$-free (consequently $K_4$-free) $t$-colorings of $K_n$, where $n\sim t^2$.  By clique substitutions, we get $C_4$-free $(t,2)$-colorings with the largest monochromatic clique having $\sim {3n\over t^2}$ vertices, thus $t$ monochromatic cliques can cover at most $m\sim {3n\over t}$ vertices. For $t=2$, Theorem \ref{fourover5} gives a sharp bound (${4n\over 5}$) for a strong cover. However, for $t=3$, we do not know the asymptotic, our best example is ${3n\over 4}$, based on the coloring of $K_8$ described in (vii) below.

\end{itemize}


\subsection{Ramsey type problems}

\begin{figure}[h!]
\begin{center}
  \includegraphics[width=0.5\linewidth]{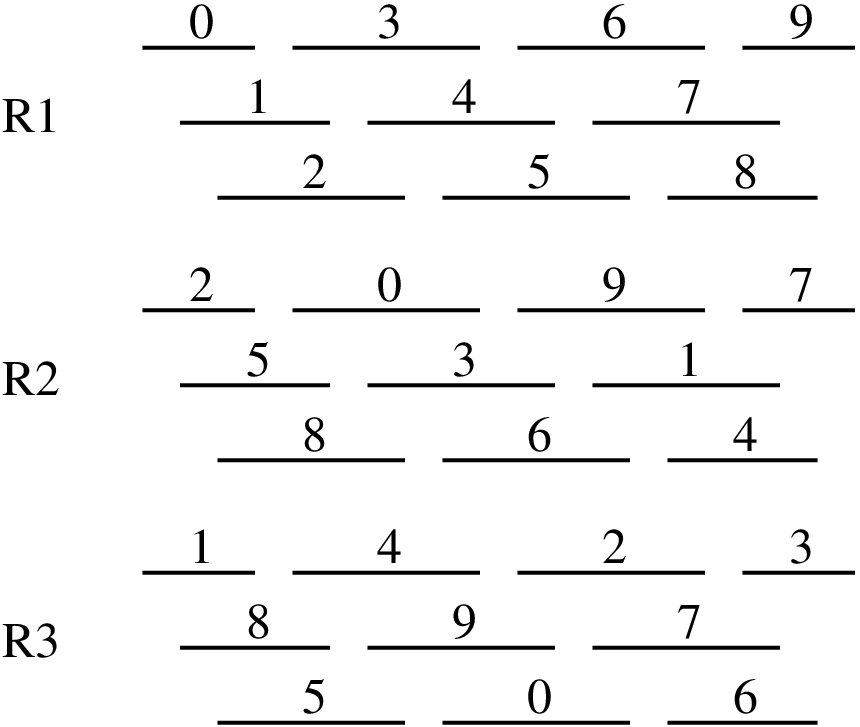}
  \caption{\label{abrach3}Covering $K_{10}$ with $3$ interval graphs. }
\end{center}
\end{figure}

\begin{itemize}
\item (iv)  In every chordal $(2,2)$-coloring of $K_n$, there is a monochromatic clique $K_m$ with $m\ge \lceil {n\over 2}\rceil$. This follows immediately from Theorem \ref{gyl2}. Equality is possible for every $n\ge 4$ with interval colorings by substituting vertices into the $(2,2)$-coloring of $K_4$, where both colors form a path with three edges.

\item (v) In every interval $(3,2)$-coloring of $K_n$, there is a monochromatic clique $K_m$ with $m\ge {n\over 4}$. This follows from the result that pairwise intersecting $3$-intervals can be pierced by four points (Theorem 4 in \cite{GYL}). On the other hand,  Figure~\ref{abrach3} shows an interval $(3,2)$-coloring of $K_{10}$ with no monochromatic $K_4$. One can blow-up each vertex by Proposition~\ref{subst} to get an interval $(3,2)$-coloring of $K_n$ where any monochromatic clique can cover at most $\lceil\frac{3n}{10}\rceil+1$ ($3n/10$ if $n=0\pmod {10}$) of the vertices.

\item (vi) For $t\ge 3$, in every chordal $(t,2)$-coloring of $K_n$, there is a monochromatic clique $K_m$ with $m\ge \alpha(t^{-1})n$. This follows from the result of \cite{AK},\cite{GYHS} cited above.  There exists an interval $t$-coloring of $K_n$ with no monochromatic clique larger than $\lceil {n\over t}\rceil+1$: partition $V(K_n)$ into $t$ parts $S_1,\dots,S_t$ evenly. For every $i$ the edges within $S_i$ and for $j>i$ the edges in the complete bipartite graph  $[S_i,S_j]$ are colored with color $i$.

\item (vii) In every $C_4$-free $(2,2)$-coloring of $K_n$, there is  a monochromatic clique $K_m$ with $m\ge \lceil {2n\over 5}\rceil$. Equality is possible for every $n\ge 5$. This follows immediately from Theorem \ref{fourover5}.

    For $C_4$-free $(3,2)$-colorings, we have no asymptotic result.  Our best lower bound is the general one from (viii). Our upper bound is $\lceil {n\over 4}\rceil+1$ ($n/4$ if $n=0\pmod 8$) obtained by substitutions into the following $3$-coloring of $K_8$ with vertex set $[8]$. Color $1$ is the $C_7$ defined by the cyclic order $(1234567)$ and the edges $(4,8),(7,8)$;  color $2$ is the $C_7$ defined by the cyclic order $(1835746)$ and the edges $(2,5),(2,6)$; color $3$ is the $C_8$ defined by the cyclic order $(14273685)$ and the diagonals  $(1,3),(2,8)$.

\item (viii)  In every $C_4$-free $(t,2)$-coloring of $K_n$, there is a monochromatic clique $K_m$ with $m\ge  \alpha^2(t^{-1})n$. This follows from the result of \cite{HO} cited above. On the other hand, as shown in (iii), there exists $C_4$-free $(t,2)$-colorings with no larger monochromatic cliques than $\sim{3n\over t^2}$ .
\end{itemize}

\end{document}